\documentstyle[12pt,bezier]{article}
\input{emlines.sty}
\input{amssym.def}	%using AmSTeX fonts, for instance:
\input{amssym}		%in $\Bbb N, N$ first character fine, second usual

\newcommand{\qed}{$\;\;\;\Box$}
\newenvironment{proof}{\par\smallbreak{\sl Proof.~}}
{\unskip\nobreak\hfill \qed \par\medbreak}

\newcommand{\D}{{\cal D}}

\newcommand{\N}{{\Bbb N}}
\newcommand{\R}{{\Bbb R}}

\newcommand{\WF}{\mathop{\mathrm{WF}}}

\renewcommand{\d}{\partial}

\newtheorem{thm}{Theorem}
\newtheorem{lemma}[thm]{Lemma}
\newtheorem{prop}[thm]{Proposition}
\newtheorem{defn}[thm]{Definition}

\newcommand{\supp}{\mathop{\mathrm{supp}}}
\newcommand{\sing}{\mathop{\mathrm{sing}}}

\title{A Distributional Solution to a
Hyperbolic Problem Arising in Population Dynamics
}

\newcounter{thesame}
\author{
Irina Kmit
\thanks{
This work was done
while visiting the Institut f{\"u}r Mathematik,
Universit{\"a}t Wien,
supported by an {\"O}AD grant.}
}

\date{}

\begin{document}

\maketitle

\begin{abstract}
We consider a generalization of the Lotka-McKendrick problem describing the 
dynamics of an age-structured population with time-dependent vital rates.
 The generalization consists in allowing the initial and the boundary 
conditions to be derivatives of the Dirac measure.  We construct
a unique $\D'$-solution in the framework of intrinsic multiplication of
distributions. We also investigate the regularity of this solution.
\end{abstract}

\emph{Key words:} population dynamics, hyperbolic equation, integral condition,
singular data, distributional solution

\emph{Mathematics Subject Classification: 35L50, 35B65, 35Q80, 58J47}

\section{Introduction}\label{sec:intr}

We consider a non-classical  hyperbolic 
problem with integral boundary condition 
\begin{eqnarray}
(\partial_t  + \partial_x) u &=& p(x,t)u + g(x,t), \qquad (x,t)\in \Pi
  \label{eq:1}\\
u|_{t=0} &=& a(x), \qquad x\in [0,L) \label{eq:2} \\
u|_{x=0} &=&  c(t)\int\limits_{0}^{L}b(x)u\,dx,
  \qquad  t\in[0,\infty) \label{eq:3} \, ,
\end{eqnarray}
where
$$
\Pi=\{(x,t)\in\R^2\,|\,0<x<L,t>0\}.
$$
From the point of view of applications, (\ref{eq:1})--(\ref{eq:3}) describes
the dynamics of the age-structured population
 (see i.e.~\cite{bu-ia,cushing,metz,Sanchez,webb}).
There $u$ denotes the distribution of individuals having age $x>0$ at time $t>0$,
$a(x)$ is the initial distribution, $-p(x,t)$  denotes the mortality rate,
$b(x)$ denotes the age-dependent fertility rate, $c(t)$ is the specific
fertility rate of females, $g(x,t)$ is the distribution
of migrants, $L$ is the maximum age attained by individuals.
Furthermore, $b(x)=0$ on $[0,L]\setminus[L_1,L_2]$, where $[L_1,L_2]\subset
[0,L]$ is the fertility period of females.
The evolution of $u$ without diffusion is governed
by (\ref{eq:1})--(\ref{eq:3}). The system (\ref{eq:1})--(\ref{eq:3})
is a continuous model of a discrete structure. As
in many problems of such a kind, it is natural to consider singular initial
and boundary data.  We focus on the case when these  data
have  singular support in finitely many points, i.e.
$$
a(x)=a_r(x)+\sum\limits_{i=1}^md_{1i}\delta^{(m_i)}(x-x_i)\quad
\mbox{for\,\,some}\quad d_{1i}\in\R,
m_i\in\N_0, x_i\in(0,L),
$$
\begin{equation}\label{eq:0}
b(x)=b_r(x)+\sum\limits_{k=1}^sd_{2i}\delta^{(n_k)}(x-x_k)\quad \mbox{for\,\,some}\quad
d_{2i}\in\R,
n_k\in\N_0, x_k\in(0,L),
\end{equation}
$$
c(t)=c_r(t)+\sum\limits_{j=1}^qd_{3i}
\delta^{(l_j)}(t-t_j)\quad \mbox{for\,\,some}\quad d_{3i}\in\R,
l_j\in\N_0, t_j\in(0,\infty).
$$
The data of the Dirac measure type enable us to model the point-concentration of various
demographic parameters. 

The problem under consideration is of interest from both biological and mathematical
points of view.

\paragraph{Biological motivation}
A basic model describing the evolution of an age-structured population is given by the
Lotka-McKendrick system:
\begin{eqnarray}
\label{eq:M}
(\partial_t  + \partial_x)u &=& -p(x)u
  \\
u|_{t=0} &=& u_0(x)\\
%\displaystyle
u|_{x=0} &=&  \int\limits_{0}^{L}b(x)u\,dx.
 \end{eqnarray}
The differential equation describes the aging of the population and the output due
to deaths. The integral $\int_{\alpha_1}^{\alpha_2}u(x,t)\,dx$ gives the number
of individuals at time $t$ having age $x$ in the range $\alpha_1\le x\le\alpha_2$. Thus,
the third equation is responsible for newborns, entering the population at age zero.

A biological generalization of~(\ref{eq:M}) to~(\ref{eq:1})--(\ref{eq:3}) consists in
the allowing the  fertility and mortality rates to depend on $t$ (see e.g.~\cite{inaba1,kim,lopez}).
In reality the vital rates are never time-homogeneous and adapt to the changing
social and technological environment. Introducing  the $\delta$-distributional data in~(\ref{eq:2}) and~(\ref{eq:3})
also has a biological meaning (see~\cite{metz}).

In demography, $c(t)$ is the total fertility rate of the population at time $t$, in other words, the
average number of childbirths per female during her reproductive period.
On one side, the results presented in the paper could shed a new light on the so-called
$c$-control problems when one wants to control the population only through changing $c(t)$.
Chinese scientists used discrete models to provide mathematical background for the
unicity child policy ($c$-control problem) in the People's Republic of China~\cite{song,songyu,yu}.
Continuous models in the context of the $c$--control problem were considered in~\cite{inaba}.
In contrast to the aforementioned papers, the presence of strongly singular data 
in~(\ref{eq:2}) and~(\ref{eq:3}) allows one to combine the continuity of the model with
the discreteness of the real evolutionary process. Occurrence
of strong singularities in $c(x)$ can be motivated by synchronized and concentrated 
reproduction of the species. This also allows one  to involve statistical data 
into (\ref{eq:1})--(\ref{eq:3}) and perhaps makes our model competitive with discrete-time 
and discrete-age models~\cite{caswell}. 

Involving strong singularities into the model could have another interpretation:
such singularities can be produced by a linearization of nonlinear problems
with discontinuous data. Thus this opens a space for 
interesting nonlinear consequences.

\paragraph{Mathematical motivation}
We consider our paper as a further
step in the study  of
generalized solutions to initial-boundary hyperbolic problems
in two variables.

Since the singularities given on $\partial\Pi$  expand
inside $\Pi$ along
characteristic curves of
the equation (\ref{eq:1}), a
solution preserves at least
the same order of regularity as it has on   $\partial\Pi$. This causes
multiplication of distributions under the integral sign in
(\ref{eq:3}). In spite of this complication, we
find distributional solutions of
(\ref{eq:1})--(\ref{eq:3}).
In parallel, we study propagation, interaction and  creation
of new singularities for the problem (\ref{eq:1})--(\ref{eq:3}).

Initial-boundary semilinear hyperbolic problems with distributional data
were studied, among others, in~\cite{13,Kmit,KmitHorm}.   There also appears
a complication with  multiplication of
distributions that
is caused by nonlinear right-hand sides
of the differential equations and also  by boundary conditions that are
nonlinear (with bounded nonlinearity) in~\cite{13},
nonseparable in~\cite{KmitHorm}, and integral in~\cite{Kmit}.
To overcome this complication, the authors
use the framework of  {\it delta waves} (see \cite{RauchReed87}).
In other words, they find solutions by regularizing all singular
data, solving the regularized system and then passing
the obtained sequential solution
to a weak limit.

Boundary and initial-boundary value problems
 for a linear second order hyperbolic equation~\cite{89} and the general
 strictly hyperbolic systems in the Leray-Volevich sense~\cite{91}
 are studied in a complete scale of {\it Sobolev type spaces} depending on
 parameters $s$ and $\tau$, where $s$ characterizes
 the smoothness of a solution in all variables and $\tau$
 characterizes additional smoothness in the tangential variables.
 Sobolev-type a priori estimates are obtained and, based on them, the
 existence and uniqueness results in Sobolev spaces are proved.

 In contrast to the aforementioned papers
% to~\cite{13,KmitHorm}
we here treat {\it integral} boundary conditions and show that the problem
(\ref{eq:1})--(\ref{eq:3})
is solvable in the {\it distributional} sense.
We construct a unique distributional solution by means of
multiplication of distributions in the sense of H\"ormander~\cite{Horm}.

We show that the
boundary condition (\ref{eq:3}) causes anomalous singularities
at the time
when singular characteristics and
vertical singular lines arising from
the data of (\ref{eq:3}) intersect.
In the case that the singular part of $b(x)$ is a sum of derivatives of
the Dirac measure,
 the solution becomes
more singular.
In the case that the initial
and the boundary data are
Dirac measures, the solution preserves the same order of regularity.
Similar phenomenon was shown in~\cite{12}  for a semilinear
hyperbolic Cauchy   problem with strongly singular initial data, where interaction 
of singularities was caused by the nonlinearity of the
equations. Anomalous singularities
were considered also in~\cite{RauchReed81} and~\cite{Ober86}, where
propagation of singularities for, respectively, initial and initial-boundary
semilinear hyperbolic problems were studied. There was proved
that, if the initial data are, at worst, jump discontinuities, then
the singularities at the common point
of singular
characteristics of the differential equations  are weaker. Furthermore,
if boundary data are regular enough, then
reflected singularities cannot be stronger than the corresponding incoming
singularities.
It turns out~\cite{Elt,LavLyu} that in some cases of nonseparable
boundary conditions the solution becomes more regular in time,
namely for $C^1$-initial data it becomes $k$-times
continuously differentiable for any desired $k\in \N_0$ in a finite time.

{\bf Organization of the paper}
Section 2 contains some
basic facts from the theory of distributions.
In Section 3 we describe our problem
in detail and state our result. Sections 4--9
present successive steps of construction of a distributional
solution to the problem. In particular, the integral boundary
condition is treated in Section~5.
 In parallel we analyze the regularity of
the solution.
 The uniqueness is proved
in Section 10.
% (\ref{eq:1})--(\ref{eq:3}).

\section{Background}
For convenience of the reader we here recall the relevant material
from~\cite{Joshi,Horm1,Horm,Shilow} without proofs. Throughout the 
paper we will denote by $\langle\cdot,\cdot\rangle: \D'\times\D\to\R$
the dual pairing on the space $\D$ of $C^\infty$-functions having 
compact support.
\begin{defn}
(\cite{Horm1}, 2.5) A distribution $u\in\D'(\R^2)$ is {\em microlocally smooth}
at $(x,t,\xi,\eta)$  $((\xi,\eta)\ne 0)$ if
the following condition holds: If  $u$ is localized
about $(x,t)$ by $\varphi\in \D(\R^2)$ with
$\varphi\equiv 1$ in a
neighborhood of $(x,t),$ then the Fourier transform of $\varphi u$
is rapidly decreasing in an open cone about $(\xi,\eta).$
The {\em wave front set} of $u$, $\WF(u)$, is the complement in $\R^4$
of the set of microlocally smooth points.
\end{defn}

\begin{prop}
(\cite{Horm}, 8.1.5)
Let $u\in\D'(\R^2)$ and $P(x,D)$ be a linear differential operator with smooth
coefficients. Then
$$
\WF(Pu)\subset \WF(u).
$$
\end{prop}

\begin{defn}
(\cite{Horm}, 6.1.2)
Let $X,Y\subset\R^2$ be open sets and $u\in\D'(Y)$. Let
$f:X\to Y$ be a smooth invertible map such that its derivative
is surjective. Then the {\em pullback of $u$ by $f$}, $f^*u$, is a unique
continuous linear map: $\D'(Y)\to\D'(X)$ such that for all
$\varphi\in\D(Y)$
$$
\left\langle f^*u,\varphi\right\rangle=\left\langle u,|J(f^{-1})|(\varphi\circ f^{-1})\right\rangle,
$$
where $J(f^{-1})$ is the Jacobian matrix of $f^{-1}$.
\end{defn}

\begin{thm}  (\cite{Horm}, 8.2.7)
Let $X$ be a manifold and $Y$ a submanifold with normal bundle denoted
by $N(Y)$. For every distribution $u$ in $X$ with $\WF(u)$ disjoint from
$N(Y)$, the restriction $u|_Y$ of $u$ to $Y$ is a well-defined distribution on
$Y$ that is the pullback by the inclusion $Y\hookrightarrow X.$
\end{thm}

\begin{thm}(\cite{Horm}, 5.1.1)
For any distributions $u\in\D'(X_1)$ and $v\in\D'(X_2)$
there exists  a unique
distribution $w\in\D'(X_1\times X_2)$ such that
$$
\left\langle
w,\varphi_1\otimes\varphi_2\right\rangle=\left\langle u,\varphi_1\right\rangle\left\langle v,\varphi_2\right\rangle,\quad
\varphi_i\in \D(X_i)
$$
and
$$
\left\langle
w,\varphi\right\rangle=\left\langle u,\left\langle v,\varphi(x_1,x_2)\right\rangle\right\rangle
=\left\langle v,\left\langle u,\varphi(x_1,x_2)\right\rangle\right\rangle,\quad
\varphi\in \D(X_1\times X_2).
$$
Here $u$ acts on $\varphi(x_1,x_2)$
as on a function of $x_1$ and  $v$ acts on $\varphi(x_1,x_2)$
as on a function of $x_2$.
\end{thm}
The
distribution $w$ as in Theorem 5 is called the {\em tensor product} of $u$
and $v$, and denoted by
$w=u\otimes v$.

\begin{thm}(\cite{Joshi}, 11.2.2)
Let $X,Y$ be open sets in $\R^2$ and let $f:X\to Y$ be a diffeomorphism.
If $u\in\D'(Y)$, then $f^*u$, the pull-back of $u$, is well defined,
and we have
$$
\WF(f^*(u))=\{(x,df_x^t\eta):(f(x),\eta)\in\WF(u)\}.
$$
\end{thm}

\begin{thm}(\cite{Horm}, 8.2.10)\label{thm:times}
If $v,w\in\D'(X)$, then the product
$v\cdot w$
is well defined as the pullback of the tensor product $v\otimes w$ by the
diagonal map   $\delta:\R\to \R\times \R$ unless
$(x,t,\xi,\eta)\in \WF(v)$ and $(x,t,-\xi,-\eta)\in \WF(w)$
for some $(x,t,\xi,\eta)$.
\end{thm}

\begin{thm}(\cite{Shilow}, 8.6)\label{thm:shilow}
If a distribution $u$ is identically equal to 0 on each of the
domains $G_i$, $i\ge 1$, then
$u$ is identically equal to 0 on $G=\bigcup\limits_{i\ge 1}G_i$.
\end{thm}

\section{Statement of the results}\label{sec:multipl}
For simplicity of technicalities we assume that both
the initial and the boundary data have singular
supports  at a single point and are the Dirac measures or
 derivatives of the Dirac measure.
This causes no loss of generality for the problem if
the singular parts of the initial and
the boundary data are finite sums of the Dirac  measures
and derivatives thereof, i.e. they are of the form~(\ref{eq:0}).
Specifically, we consider the following system
\begin{eqnarray}
(\partial_t  + \partial_x) u &=& p(x,t)u+g(x,t),
  \label{eq:51}\\
u|_{t=0} &=& a_r(x)+\delta^{(m)}(x-x_1^*), \qquad x\in [0,L) \label{eq:52} \\
u|_{x=0} &=& (c_r(t)+\delta^{(j)}(t-t_1))\int\limits_{0}^{L}
(b_r(x)+\delta^{(n)}(x-x_1))u\,dx,
   \label{eq:53} \, \\ &&t\in[0,\infty),\nonumber
\end{eqnarray}
where $x_1>0,x_1^*>0,t_1>0$, and
$m,j,n\in\N_0$.  Without loss of generality we can assume that  $x_1^*<x_1$.

We impose the following conditions:

\vskip1.0mm
{\it Assumption 1.}
$a_r^{(i)}(0)=0,  c_r^{(i)}(0)=0$ for all $i\in\N_0.$

\vskip1.0mm
{\it Assumption 2.}
$b_r^{(i)}(L)=0$ for all    $i\in\N_0$ and there exists $\varepsilon>0$
such that  $b_r(x)=0$ for
$x\in[0,\varepsilon]$.

\vskip1.0mm
{\it Assumption 3.}
The functions  $p$ and $g$ are  smooth in
$\R^2$, $a_r$ is smooth on $[0,L)$,  $b_r$ is  smooth on $[0,L]$,
and $c_r$ is  smooth
on $[0,\infty)$.
 \vskip1.0mm

Note that Assumption 1 ensures the arbitrary order
compatibility between (\ref{eq:52}) and (\ref{eq:53}). Assumption~2
 not particularly restrictive from the practical point of
view, since $[0,L]$ covers the fertility period of females.

All characteristics of the differential equation (\ref{eq:1}) as solutions to
the following initial problem for ordinary differential equation:
$$
\frac{dx}{dt}=1,\quad x(t_0)=x_0,
$$
where $(x_0,t_0)\in\R^2$, are given by the formula $x=t+x_0-t_0$.

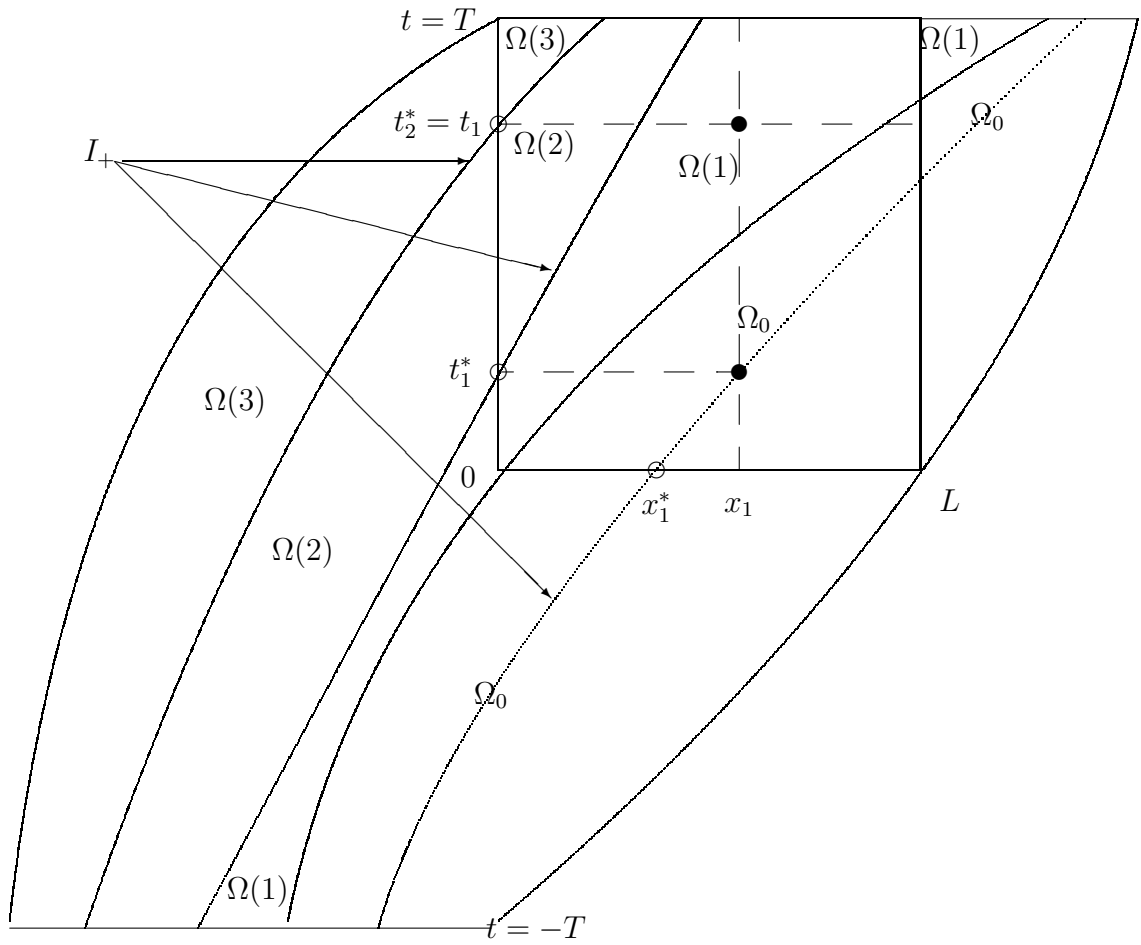
\break

\begin{figure}[t]
\centerline{
\unitlength=1.00mm
\special{em:linewidth 0.4pt}
\linethickness{0.4pt}
\begin{picture}(155.00,125.00)
\put(70.00,65.00){\framebox(56.00,60.00)[cc]{}}
\bezier{616}(155.00,125.00)(141.00,66.00)(70.00,5.00)
\bezier{684}(42.00,5.00)(55.00,72.00)(143.00,125.00)
\bezier{616}(70.00,125.00)(15.00,96.00)(5.00,5.00)
\emline{5.00}{4.00}{1}{69.00}{4.00}{2}
\emline{126.00}{125.00}{3}{155.00}{125.00}{4}
\bezier{300}(148.00,125.00)(66.00,47.00)(54.00,4.00)
\emline{102.00}{65.00}{5}{102.00}{68.00}{6}
\emline{102.00}{71.00}{7}{102.00}{75.00}{8}
\emline{102.00}{79.00}{9}{102.00}{83.00}{10}
\emline{102.00}{87.00}{11}{102.00}{92.00}{12}
\emline{102.00}{96.00}{13}{102.00}{100.00}{14}
\emline{102.00}{104.00}{15}{102.00}{109.00}{16}
\emline{102.00}{113.00}{17}{102.00}{118.00}{18}
\emline{102.00}{122.00}{19}{102.00}{125.00}{20}
\emline{102.00}{78.00}{21}{99.00}{78.00}{22}
\emline{99.00}{78.00}{23}{99.00}{78.00}{24}
\emline{96.00}{78.00}{25}{92.00}{78.00}{26}
\emline{88.00}{78.00}{27}{84.00}{78.00}{28}
\emline{80.00}{78.00}{29}{76.00}{78.00}{30}
\emline{72.00}{78.00}{31}{70.00}{78.00}{32}
\bezier{600}(97.00,125.00)(81.00,99.00)(30.00,4.00)
\emline{70.00}{111.00}{33}{73.00}{111.00}{34}
\emline{76.00}{111.00}{35}{80.00}{111.00}{36}
\emline{84.00}{111.00}{37}{88.00}{111.00}{38}
\emline{91.00}{111.00}{39}{95.00}{111.00}{40}
\emline{98.00}{111.00}{41}{102.00}{111.00}{42}
\emline{105.00}{111.00}{43}{109.00}{111.00}{44}
\emline{113.00}{111.00}{45}{117.00}{111.00}{46}
\emline{121.00}{111.00}{47}{125.00}{111.00}{48}
\bezier{600}(84.00,125.00)(45.00,91.00)(15.00,4.00)
\put(75.00,4.00){\makebox(0,0)[cc]{$t=-T$}}
\put(62.00,125.00){\makebox(0,0)[cc]{$t=T$}}
\put(130.00,61.00){\makebox(0,0)[cc]{$L$}}
\put(66.00,64.00){\makebox(0,0)[cc]{$0$}}
\put(91.00,60.00){\makebox(0,0)[cc]{$x_1^*$}}
\put(102.00,60.00){\makebox(0,0)[cc]{$x_1$}}
\put(65.00,78.00){\makebox(0,0)[cc]{$t_1^*$}}
\put(62.00,111.00){\makebox(0,0)[cc]{$t_2^*=t_1$}}
\put(91.00,65.00){\circle{2.00}}
\put(70.00,78.00){\circle{2.00}}
\put(70.00,111.00){\circle{2.00}}
\put(102.00,78.00){\circle*{2.00}}
\put(102.00,111.00){\circle*{2.00}}
\put(104.00,85.00){\makebox(0,0)[cc]{$\Omega_0$}}
\put(135.00,112.00){\makebox(0,0)[cc]{$\Omega_0$}}
\put(69.00,35.00){\makebox(0,0)[cc]{$\Omega_0$}}
\put(98.00,105.00){\makebox(0,0)[cc]{$\Omega(1)$}}
\put(130.00,122.00){\makebox(0,0)[cc]{$\Omega(1)$}}
\put(38.00,9.00){\makebox(0,0)[cc]{$\Omega(1)$}}
\put(72.00,108.00){\makebox(0,0)[lc]{$\Omega(2)$}}
\put(44.00,54.00){\makebox(0,0)[cc]{$\Omega(2)$}}
\put(75.00,122.00){\makebox(0,0)[cc]{$\Omega(3)$}}
\put(35.00,74.00){\makebox(0,0)[cc]{$\Omega(3)$}}
\put(19.00,106.00){\vector(1,-1){58.00}}
\put(19.00,106.00){\vector(4,-1){58.00}}
\put(20.00,106.00){\vector(1,0){46.00}}
\put(17.00,107.00){\makebox(0,0)[cc]{$I_+$}}
\end{picture}
}
\end{figure}

\mbox{}

\break

\begin{defn}
Let $I_+=\bigcup_{n\ge 0}I_+[n]$,
where $I_+[n]$ are
subsets of $\R^2$ defined
by induction as follows.

\begin{itemize}
\item
 $I_+[0]$ is the union of the characteristics $x=t+x_1^*$
and $x=t-t_1$.

\item
Let $n\ge 1$. If $I_+[n-1]$ includes the characteristic
 $x=t+x_1-\tilde t$, then
$I_+[n]$ includes  the characteristic $x=t-\tilde t$.
\end{itemize}
\end{defn}
For characteristics contributing into
$I_+$
denote their intersection points   with the positive semiaxis $x=0$
by $t_1^*,t_2^*,\dots$.
We assume that
$t_j^*<t_{j+1}^*$ for $j\ge 1$.
 The union of all singular characteristics
of the initial problem, as it will be shown,
is included into the set $I_+$. In fact, we will show that
$\sing\supp u\subset I_+$.

\vskip1.0mm
{\it Assumption 4.}
$
x_1-t_1\ne x_1^*, \quad
t_1-x_1\ne t_s^*$
for all $t_s^*<t_1.$
 \vskip1.0mm
This assumption excludes the
situation when three different singularities
intersect at the same point. Without this assumption
the distributional solution does not exist, because there appears
 multiplication of the Dirac measure onto itself.

Our goal is, using distributional multiplication, to obtain
distributional solution to (\ref{eq:51})--(\ref{eq:53}).
We use the notion of the so-called ''WF favorable'' product
which is due to L.~H\"ormander \cite{Horm} and is in the second level
 of M.~Ober\-gug\-genberger's
hierarchy of intrinsic distributional products \cite[p.~69]{11}.

We actually obtain distributional solution  in a domain
$\Omega\subset\R^2$ that is the  domain of influence (or determinacy)
of the problem (\ref{eq:51})--(\ref{eq:53}).
 Clearly, $\Omega$ is the union of all characteristics
 $x=t+x_0-t_0$ passing through those points $(x_0,t_0)$
 on the boundary of $\Pi$ where the conditions
 (\ref{eq:52}) and (\ref{eq:53}) are given, i.e. through points
 $(x_0,t_0)\in\left([0,L)\times\{0\}\right)\cup\left(\{0\}\times[0,\infty)\right)$. In other words,
 $$
 \Omega=\{(x,t)\in\R^2\,|\,x<t+L\}.
 $$

 \begin{defn}\label{defn:Omega}
 A distribution $u$ is called a $\D'(\Omega)$-solution to
 the problem (\ref{eq:51})--(\ref{eq:53}) if the following conditions are
met.\\
 1. The equation (\ref{eq:51}) is satisfied in $\D'(\Omega)$: for every
 $\varphi\in\D(\Omega)$
 $$
\left\langle(\d_t+\d_x-p(x,t))u,\varphi\right\rangle=\left\langle  g(x,t),\varphi\right\rangle.
 $$
 2. $u$ is restrictable to $[0,L)\times\{0\}$
 in the sense of H\"ormander (see Theorem 4) and
 $u|_{t=0} = a_r(x)+\delta^{(m)}(x-x_1^*), \quad x\in [0,L)$.\\
 3. The product of $(b_r(x)+
\delta^{(n)}(x-x_1))\otimes 1(t)$ and $u(x,t)$
 exists in $\D'(\Pi)$ in the sense of H\"ormander (see Theorem 7).\\
 4. $\int_{0}^{L}\left[\left(b_r(x)+\delta^{(n)}(x-x_1)\right)\otimes 1(t)\right]u\,dx$  is a distribution
$v\in\D'(\R_+)$ defined by
$$
\left\langle v,\psi(t)\right\rangle=\left\langle[(b_r(x)+
\delta^{(n)}(x-x_1))\otimes 1(t)]u,1(x)\otimes
\psi(t)\right\rangle,\quad \psi(t)\in\D(\R_+),
$$
where $b_r(x)=0$, $x\notin[0,L]$.\\
5. $v$ is a smooth function in $t_1$.\\
6. $u$ is restrictable to  $\{0\}\times[0,\infty)$
 in the sense of H\"ormander (see Theorem 4) and
 $u|_{x=0} = (c_r(t)+\delta^{(j)}(t-t_1))v$,\quad $t\in[0,\infty)$.\\
7. $\sing\supp u\subset\Omega\setminus\{(x,t)\,|\,x=t\}$.
 \end{defn}

 Our next objective is to define the solution concept for
 (\ref{eq:51})--(\ref{eq:53}) on $\Pi$. It is not so obvious how we should
 define the restriction of $u\in\D'(\Pi)$ to the boundary of
 $\Pi$ so that the initial and the boundary conditions are meaningful.
In this respect let us make the following observation.

 Note that $\overline{\Pi}\setminus\{(L,0)\}\subset\Omega$. Let
 $\Omega_0\subset\Omega$ be a domain  such that
 $\overline{\Pi}\setminus\{(L,0)\}\subset\Omega_0$ and $u$ be a
$\D'(\Omega)$-solution  to the problem (\ref{eq:51})--(\ref{eq:53})
in the sense of Definition~\ref{defn:Omega}. Then $u$ restricted
to $\Omega_0$ is a $\D'(\Omega_0)$-solution
to the problem (\ref{eq:51})--(\ref{eq:53})
in the sense of the same definition. This suggests the following definition.

\begin{defn}\label{defn:Pi}
Let $u$ be a $\D'(\Omega)$-solution to the problem
(\ref{eq:51})--(\ref{eq:53})
in the sense of Definition~\ref{defn:Omega}.
Then $u$ restricted to $\Pi$ is called
a $\D'(\Pi)$-solution to the problem~(\ref{eq:51})--(\ref{eq:53}).
\end{defn}

Set
$$
\Omega_+=\left\{(x,t)\in\Omega\,
|\,x>0,t>0\right\}.
$$
We are now prepared to state the existence result.
\begin{thm}\label{thm:exist}
1. Let Assumptions 1--4 hold. Then there exists a
$\D'(\Omega)$-solution~$u$ to the problem
(\ref{eq:51})--(\ref{eq:53}) in the sense of Definition~\ref{defn:Omega}
such that:
\begin{equation}\label{eq:r}
\parbox{11cm}{the restriction of $u$ to
 any domain $\Omega_+'\supset\Omega_+$
such that any characteristic of
(\ref{eq:51}) intersects $\d\Omega_+'$
at a single point does not depend
on the values of the functions
$p$ and $g$
on
$\Omega\setminus\Omega_+'$.}
\end{equation}
\\
2. Let Assumptions 1--4 hold. Then there exists a    $\D'(\Pi)$-solution to the problem
(\ref{eq:51})--(\ref{eq:53}) in the sense of Definition~\ref{defn:Pi}.
\end{thm}

Given a domain $G$, set
$$
\D_+'(G)=\left\{u\in\D'(G)\,|\,
u=0\,\,\mbox{for}\,\,
(x<0\,\,\mbox{or}\,\,t<0)\right\}.
$$
Write
$$
\tilde\lambda(x,t)=
\cases{1
&if
$(x,t)\in\overline{\Omega_+}$,
 \cr 0
&if
$(x,t)\in
\Omega\setminus\overline{\Omega_+}$, \cr}
$$
$$
\tilde p(x,t)=
\cases{p
&if
$(x,t)\in\overline{\Omega_+}$,
 \cr 0
&if
$(x,t)\in
\Omega\setminus\overline{\Omega_+}$. \cr}
$$
Similarly to  $p$ we define a modification of  $g$ and denote
it by $\tilde g$.

\begin{defn}\label{defn:D+}
 $u\in\D_+'(\Omega)$ is called a $\D_+'(\Omega)$-solution to the problem
(\ref{eq:51})--(\ref{eq:53}) if the following conditions are met.\\
1. Items 3--5 of Definition~\ref{defn:Omega} hold.\\
2. The equation (\ref{eq:51}) is satisfied in $\D_+'(\Omega)$: for every
 $\varphi\in\D(\Omega)$
\begin{eqnarray*}
\lefteqn{
 \left\langle(\d_t+\tilde\lambda(x,t)\d_x-\tilde p(x,t))u,\varphi\right\rangle=\left\langle\tilde g(x,t),\varphi\right\rangle}\\
&&
+\left\langle(a_r(x)+\delta^{(m)}(x-x_1^*))\otimes\delta(t)+
\delta(x)\otimes[(c_r(t)+
\delta^{(j)}(t-t_1))v],\varphi\right\rangle,
 \end{eqnarray*}
where $a_r(x)=0$ if $x<0$ and $v(t)=0$ if $t\le 0$.
\\
3. $\sing\supp u\setminus\d\Omega_+\subset\Omega_+\setminus
\{(x,t)\,|\,x=t\}$.
\end{defn}

\begin{prop}\label{prop:p}
Let $u$ be a $\D'(\Omega)$-solution to the problem
(\ref{eq:51})--(\ref{eq:53})
in the sense of Definition~\ref{defn:Omega}
that satisfies
(\ref{eq:r}).
Then there exists a
$\D_+'(\Omega)$-solution $\tilde u$ to the problem
(\ref{eq:51})--(\ref{eq:53}) in the sense of Definition~\ref{defn:D+}
such that
$$
u=\tilde u \,\,\mbox{in}\,\, \D'(\Omega_+).
$$
\end{prop}
This proposition is a straightforward consequence of
Definitions~\ref{defn:Omega} and~\ref{defn:D+}.
Since $\Pi\subset\Omega_+$, it makes sense  to state the uniqueness result in
$\D_+'(\Omega)$.
Write
\begin{equation}\label{eq:S}
\begin{array}{c}
\displaystyle
S(x,t)
=\exp\left\{\int\limits_{\theta(x,t)}^tp(\tau;x,t,\tau)\,d\tau
\right\},
\end{array}
\end{equation}
where $\theta(x,t)=(t-x)H(t-x)$, $H(z)$ is the Heaviside function.
We write $\hat S$ for the function $S$
given by (\ref{eq:S}), where $p$ is replaced
by $-p$.

{\it Assumption 5.}
For every $T_0>0$ there exists $T>T_0$ such that
$\hat S(x,T)\ne 0$ for all $x$
such that $(x,T)\in\Omega_+$.
\vskip1.0mm

\begin{thm}\label{thm:uniq}
1. Let Assumptions 1--5 hold. Then
a $\D_+'(\Omega)$-solution
to the problem
(\ref{eq:51})--(\ref{eq:53})
is unique.\\
2.
 Let Assumptions 1--5 hold. Then
a $\D'(\Pi)$-solution
to the problem
(\ref{eq:51})--(\ref{eq:53})
is unique.
\end{thm}

From the construction of a $\D'(\Omega)$-solution presented in
the proof of Theorem~\ref{thm:exist} we will see that in general there appear
new singularities stronger than the initial singularities. In other words,
the singular
order (cf.~\cite[\S 13]{Shilow}) of the distributional solution
grows in time. We state this result in the following theorem.

\begin{thm}\label{thm:order}
1. Let $u$ be the $\D'(\Pi)$-solution
to the problem
(\ref{eq:51})--(\ref{eq:53}), where $n\ge 1$ and $S(x,t)\ne 0$ for all
$(x,t)\in\Pi$.
Then for each $i\ge 1$ there exist
$j>i$ and $n'\ge 1$  such that
the singular order of $u$ is equal
to $n'$ in a
neighborhood of  $x=t-t_i^*$ and
the singular order of $u$ is equal to
$n'+n$ in a
neighborhood of $x=t-t_j^*$.\\
2. If $n=j=m=0$, then the singular order of $u$ on
$\Pi$ is equal to 1.
\end{thm}

We now start with the proof of Theorem~\ref{thm:exist} which will
take Sections 4--9.
 It is sufficient to solve the problem in the domain
$$
\Omega^T=
\{(x,t)\in\Omega\,|\,t-T<x, -T<t<T\}
$$
(see the picture) for an
arbitrary fixed $T>0.$ Observe that $\Omega^T$ is the intersection of the
strip $\R\times(-T,T)$ with the domain of determinacy of (\ref{eq:51})
with respect to the set
$([0,L)\times\{0\})\cup(\{0\}
\times[0,T))$.
Fix $T>0$. We start with a subdomain
$$
\Omega_0^T=\{(x,t)\in\Omega^T\,|\,t<x<t+L\}.
$$

\section{The solution on $\Omega_0^T$}
Observe that $\Omega_0^T$ is the intersection of the strip
$\R\times(-T,T)$ with the domain of  determinacy of the
problem (\ref{eq:51})--(\ref{eq:52}).
In the case that the initial data are functions, a unique solution to the
problem (\ref{eq:51})--(\ref{eq:52})
on $\Omega_0^T$
can be written in the form
\begin{equation}\label{eq:54}
u(x,t)=S_1(x,t)+S(x,t)a_r(x-t)+S(x,t)\delta^{(m)}(x-t-x_1^*)
\end{equation}
with the functions $S(x,t)$ given by
(\ref{eq:S}) and
\begin{eqnarray}
\label{eq:S1}
\lefteqn{
S_1(x,t)
=\exp\left\{\int\limits_{\theta(x,t)}^tp(\tau+x-t),\tau\,d\tau\right\}}\nonumber\\
&&
\times
\int\limits_{\theta(x,t)}^t
\exp\left\{-\int\limits_{\theta(x,t)}^{\tau}p(\tau_1+x-t,\tau_1)\,d\tau_1
\right\}g(\tau+x-t,\tau)\,d\tau.
\end{eqnarray}
Let $A_i(x,t)=\delta^{(i)}(x)\otimes 1(t)$ and
$B_i(x,t)=1(x)\otimes \delta^{(i)}(t)$ be  the distributions in $\R^2$ that are
derivatives of the Dirac measure $\delta^{(i)}(x)$ and
$\delta^{(i)}(t)$ supported along the $t$-axis and the $x$-axis, respectively.
They are defined
by the equalities
$$
\left\langle A_i(x,t),\varphi(x,t)\right\rangle=(-1)^i\int\varphi_x^{(i)}(0,t)\,dt,
$$
$$
\left\langle B_i(x,t),\varphi(x,t)\right\rangle=(-1)^i
\int\varphi_t^{(i)}(x,0)\,dx
$$
for all $\varphi\in\D(\R^2)$.
When $i=0$, then we have the Dirac measure supported along the respective axes.

Let $f$ be the
smooth map
$$
f:(x,t)\to (x,x-t-x_1^*).
$$
The inverse
$$
f^{-1}:(x,t)\to (x,x-t-x_1^*)
$$
is unique and maps the $x$-axis to the curve $t=x-x_1^*$ and the $t$-axis onto itself.
Moreover,
$$
  f'(x,t)=\left(
  \begin{array}{cc}
    1&0\\
    1&-1
  \end{array}
  \right).
$$
 Hence the Jacobian of $f$
$$
J(f)=|f'|=-1\ne 0
$$
and  $f^*B_m=\delta^{(m)}(x-t-x_1^*)$, the pullback of
$B_m$ by $f$ (see Definition 3), is well defined. Therefore the distribution
$\delta^{(m)}(x-t-x_1^*)$ acts on test functions
$\varphi\in\D(\R^2)$ in the following way:
$$
\left\langle\delta^{(m)}(x-t-x_1^*),\varphi(x,t)\right\rangle=\left\langle f^*B_m,\varphi(x,t)\right\rangle
=-\left\langle B_m,\varphi(x,t)\circ f^{-1}(x,t)\right\rangle
$$
$$
=(-1)^{m+1}\int\d_t^m\varphi(x,x-t-x_1^*)\Big|_{t=0}\,dx
=-\int \d_t^m\varphi(x,t)|_{t=x-x_1^*}\,dx.
$$
Hence, similarly to  $B_m$,
$f^*B_m$ is the $m$-th derivative of the
Dirac measure supported along the line $t=x-x_1^*$.

 \begin{defn}\label{defn:Omega0}
 A distribution $u$ is called a $\D'(\Omega_0^T)$-solution to
 the problem (\ref{eq:51}),
(\ref{eq:52})
if
Items 1 and 2 of Definition~\ref{defn:Omega} with $\Omega$ replaced by
$\Omega_0^T$ hold.
\end{defn}

\begin{lemma}\label{lemma:0}
$u(x,t)$ given by the formula (\ref{eq:54}) is a  $\D'(\Omega_0^T)$-solution
to the problem (\ref{eq:51})--(\ref{eq:52}).
\end{lemma}
\begin{proof}
A straightforward verification shows that the
sum of the first two summands in (\ref{eq:54}) is a smooth (and, therefore, distributional)
solution to the problem (\ref{eq:51})--(\ref{eq:52}) with
the singular part of the initial condition (\ref{eq:52})  identically equal to
0. Our  goal is now to prove that the third summand in
(\ref{eq:54}) is a distributional solution to the homogeneous equation (\ref{eq:51}) with singular initial
condition $\delta^{(m)}(x-x_1^*)$.
Indeed, for all $\varphi\in\D(\Omega_0^T)$, we have
\begin{eqnarray*}
\lefteqn{
\left\langle(\d_t+\d_x)(S\delta^{(m)}(x-t-x_1^*)),\varphi\right\rangle}\\
&&
=-\left\langle S\delta^{(m)}(x-t-x_1^*),\d_t\varphi+\d_x\varphi\right\rangle\\
&&
=-\left\langle\delta^{(m)}(x-t-x_1^*),S\d_t\varphi+S\d_x\varphi\right\rangle\\
&&
=-\left\langle\delta^{(m)}(x-t-x_1^*),\d_t(S\varphi)+\d_x(S\varphi)-\d_tS\varphi-\d_xS\varphi\right\rangle.
\end{eqnarray*}
Since $w=\delta^{(m)}(x-t-x_1^*)$ is a distribution in $x-t$,
this is a weak solution
to the equation $(\d_t+\d_x)w=0$. Note that $S\varphi\in\D(\Omega_0^T)$. Therefore
$$
\left\langle\delta^{(m)}(x-t-x_1^*),\d_t(S\varphi)+\d_x(S\varphi)\right\rangle=0.
$$
By (\ref{eq:S}),
$
\d_tS+\d_xS=pS.
$
The desired assertion is therewith proved.
%This completes the proof that the third summand in (\ref{eq:54})
%is a distributional solution to the homogeneous equation (\ref{eq:51}).

It remains to prove that $S(x,t)\delta^{(m)}(x-t-x_1^*)$ may be restricted to the initial
interval $X=[0,L)\times\{0\} $. For this purpose we use  Theorems 4 and 6.
Observe that $f$ restricted to $\Omega_0^T$ is a diffeomorphism.
We check the
condition
\begin{equation}\label{eq:57}
\WF(Sf^*B_m)\cap N(X)=\emptyset,
\end{equation}
where the normal bundle $N(X)$ to $X$
is defined by the formula
$$
N(X)=\{(x,t,\xi,\eta)\,|\,(x,t)\in X, \left\langle T_{(x,t)}(X),(\xi,\eta)\right\rangle=0\}
$$
and $T_{(x,t)}(X)$ is the space of all tangent vectors to $X$ at $(x,t)$. It is
clear that in our case
$$
N(X)=\{(x,0,0,\eta),\eta\ne 0\}.
$$
Let us now look at $\WF(Sf^*B_m)$. By Proposition 2, we have
$$
\WF(Sf^*B_m)\subset \WF(f^*B_m).
$$
Recall that by definition
$$
\WF(f^*B_m)=\{(x,t,df_x^t\cdot(\xi,\eta)): (f(x,t),\xi,\eta)\in \WF(B_m)\}.
$$
We also have
$$
\WF(B_m)\subset \WF(B_0)=\{(x,0,0,\eta),\eta\ne 0\}.
$$
It follows that  $f(x,t)$  is
equal to $(x,0)$. Therefore $(x,t)=(x,x-x_1^*)$.
Furthermore,
$$
df_x^t=\left(
  \begin{array}{cc}
    1&1\\
    0&-1
  \end{array}
  \right),
\qquad
df_x^t\cdot(0,\eta)=\left(
  \begin{array}{c}
   \eta\\-\eta
  \end{array}
  \right).
$$
As a consequence,
$$
\WF(Sf^*B_m)\subset
  \{(x,x-x_1^*,\eta,-\eta), \eta\ne 0\}.
$$
This means that
$S(x,t)\delta^{(m)}(x-t-x_1^*)$
is restrictable to $X$. Considering the distribution $\delta^{(m)}(x-t-x_1^*)$ to be
smooth in $t$ with distributional values in~$x$, the initial condition
(\ref{eq:57}) follows from (\ref{eq:54}).
This finishes the proof.
\end{proof}

We have proved
that $u$ defined by (\ref{eq:54})
satisfies Items~1 and 2 of Definition~\ref{defn:Omega}
with $\Omega$ replaced by $\Omega_0^T$.
Items~4--7 on $\Omega_0^T$ do not need any proof. Item~3 will be given by
Lemma~\ref{lemma:item4} in the next section.

\section{Multiplication of distributions under the integral 
in~(\protect\ref{eq:3})}
In the further sections we will extend
 the solution over
$$
\Omega_1^T=\{(x,t)\in\Omega^T\,|\,t-T<x<t\}.
$$
We use the fact that any $\D'(\Omega)$-solution $u$ 
to our problem 
is representable as
\begin{equation}\label{eq:u}
u(x,t)=u_0(x,t)+u_1(x,t),
\end{equation}
where $u_0=u$ in $\D'(\Omega_0^T)$,
$u_0$ is identically equal to 0 on
$\overline{\Omega_1^T}$,
$u_1=u$ in $\D'(\Omega_1^T)$,  and
$u_1$ is identically equal to 0 on
$\overline{\Omega_0^T}$.
 Indeed,  if $u$ is a solution, then it is a smooth function 
in a neighborhood of 
$\{(x,t)\,|\,x=t\}$ (see Item 7 of Definition~\ref{defn:Omega}). For an arbitrary $\varphi\in\D(\Omega^T)$
consider a  representation
$
\varphi(x,t)=\varphi_1(x,t)+\varphi_2(x,t)+\varphi_3(x,t)
$
such that $\varphi_i(x,t)\in\D(\Omega^T)$, $\supp\varphi_1\subset\Omega_0^T$,
$\supp\varphi_2\cap\sing\supp u=\emptyset$, and $\supp\varphi_3\subset\Omega_1^T$.
Hence
\begin{eqnarray*}
\lefteqn{
\left\langle u_0+u_1,\varphi\right\rangle
=\left\langle u_0,\varphi_1+\varphi_2\rangle+\langle u_1,\varphi_2+\varphi_3\right\rangle}\\
&&
=\left\langle u,\varphi_1\right\rangle+\left\langle u_0,\varphi_2\right\rangle+
\left\langle u_1,\varphi_2\right\rangle+\left\langle u,\varphi_3\right\rangle
=\left\langle u,\varphi_1+\varphi_2+\varphi_3\right\rangle=\left\langle u,\varphi\right\rangle.
\end{eqnarray*}

Using (\ref{eq:u}), we rewrite $v(t)$ (see
Item~4 of Definition~\ref{defn:Omega}) in the form:
$$
v(t)=\int\limits_{0}^Lb(x)u_0(x,t)\,dx+\int\limits_{0}^Lb(x)u_1(x,t)\,dx.
$$
In this section we compute
the integral
\begin{equation}\label{eq:I_0}
I_0(t)=\int\limits_{0}^Lb(x)u_0(x,t)\,dx,\quad 0<t<T,
\end{equation}
that will be used in the construction. We have to tackle
the multiplication of distributions
involved in the integrand. For technical reasons
we
extend $a_r(x)$ and $b_r(x)$ over all $\R$ defining them to be $0$
outside $[0,L]$.
By (\ref{eq:54}),
we rewrite (\ref{eq:I_0}) as follows
\begin{eqnarray*}
\lefteqn{
I_0(t)
=\int\limits_t^Lb_r(x)(S(x,t)a_r(x-t)+S_1(x,t))\,dx}\\
&&
+
\int\limits_{0}^L
\delta^{(n)}(x-x_1)(S(x,t)a_r(x-t)+S_1(x,t))\,dx\\
&&
+\int\limits_{0}^L b_r(x)S(x,t)\delta^{(m)}(x-t-x_1^*)\,dx\\
&&
+\int\limits_{0}^L\delta^{(n)}(x-x_1)S(x,t)\delta^{(m)}(x-t-x_1^*)\,dx.
\end{eqnarray*}
To compute the second integral we take a test function $\psi(t)\in\D(0,T)$
and consider the dual pairing (see Definition~\ref{defn:Omega}, Item 4)
\begin{eqnarray*}
\lefteqn{
\left\langle\delta^{(n)}(x-x_1)(S(x,t)a_r(x-t)+S_1(x,t)),1(x)\otimes\psi(t)\right\rangle}\\
&&
=\left\langle\delta^{(n)}(x-x_1)\otimes 1(t),(S(x,t)a_r(x-t)+S_1(x,t))\psi(t)\right\rangle\\
&&
=(-1)^n\left\langle 1(t),(S(x,t)a_r(x-t)+S_1(x,t))_{x}^{(n)}|_{x=x_1}\psi(t)\right\rangle\\
&&
=(-1)^n\left\langle (S(x,t)a_r(x-t)+S_1(x,t))_{x}^{(n)}|_{x=x_1},\psi(t)\right\rangle.
\end{eqnarray*}
Let us compute the third integral:
\begin{eqnarray*}
\lefteqn{
\left\langle S(x,t)b_r(x)
\delta^{(m)}(x-t-x_1^*),
1(x)\otimes\psi(t)\right\rangle}\\
&&
=\left\langle q^*\delta^{(m)}(x),S(x,t)b_r(x)\psi(t)\right\rangle\\
&&
=\left\langle\delta^{(m)}(x),
\left(S(x+t+x_1^*,t)b_r(x+t+x_1^*)\psi(t)\right)\circ q^{-1}\right\rangle\\
&&
=(-1)^m\left\langle 1(t),\d_x^m\left(S(x+t+x_1^*,t)b_r(x+t+x_1^*)\right)\Big|_{x=0},
\psi(t)\right\rangle\\
&&
=(-1)^m\left\langle S(x+t+x_1^*,t)b_r(x+t+x_1^*)
\Big|_{x=0},\psi(t)\right\rangle.
\end{eqnarray*}
To compute the last integral in the expression for $I_0(t)$ we need the following
fact.
\begin{lemma}\label{lemma:vw}
The product of two distributions $v=\delta^{(n)}(x-x_1)\otimes 1(t)$ and $w=\delta^{(m)}(x-t-x_1^*)$
exists in the sense of H\"ormander (see Theorem~\ref{thm:times}).
\end{lemma}

\begin{proof}
 Recall that
$$
\WF(v)=\{(x_1,t,\xi_1,0),\xi_1\ne 0\},
$$
$$
\WF(w)\subset\{x,x-x_1^*,\xi_2,\xi_2,-\xi_2\ne 0\}.
$$
Thus all conditions of Theorem 7 are true and
the lemma follows.
\end{proof}

We have proved the following lemma.

\begin{lemma}\label{lemma:item4}
A distribution
$u$ defined by (\ref{eq:54})
satisfies Item 4 of  Definition~\ref{defn:Omega} with $\Pi$ replaced
by $\Pi\cap\Omega_0^T$.
\end{lemma}

Turning back to computing the last integral in $I_0(t)$,
consider the map
$$
H:(x,t)\to(x-x_1,x-t-x_1^*)
$$
and the inverse map
$$
H^{-1}:(x,t)\to(x+x_1,x-t+x_1-x_1^*).
$$
Define $H^*A_n=\delta^{(n)}(x-x_1)\otimes 1(t)$ and
$H^*B_m=\delta^{(m)}(x-t-x_1^*)$. Let us check
that the former definition is unambiguous: For any $\varphi\in\D(\R^2)$ we have
\begin{eqnarray*}
\lefteqn{
\left\langle H^*A_n,\varphi(x,t)\right\rangle=-\left\langle A_n,\varphi(x+x_1,x-t+x_1-x_1^*)
\right\rangle}\\
&&
=-\left\langle\delta^{(n)}(x),\int \varphi(x+x_1,x-t+x_1-x_1^*))\,dt\right\rangle\\
&&
=\left\langle\delta^{(n)}(x),\int \varphi(x+x_1,\tau)\,d\tau\right\rangle\\
&&
=(-1)^n\int \varphi_x^{(n)}(x_1,\tau)\,d\tau=\left\langle\delta^{(n)}(x-x_1)\otimes 1(t),\varphi(x,t)\right\rangle.
\end{eqnarray*}
Here we used a simple change of coordinates $t\to\tau=x-t+x_1-x_1^*$.

We are now in a position to compute the product of two distributions $\delta^{(n)}(x-x_1)$ and
$\delta^{(m)}(x-t-x_1^*)$: For any $\varphi\in\D(\R^2)$ we have
\begin{eqnarray*}
\lefteqn{
\left\langle S(x,t)\delta^{(n)}(x-x_1)\delta^{(m)}(x-t-x_1^*),\varphi(x,t)\right\rangle}\\
&&
=\left\langle H^*A_nH^*B_m,S(x,t)\varphi(x,t)\right\rangle\\
&&
=\left\langle H^*(A_nB_m),S(x,t)\varphi(x,t)\right\rangle\\
&&
=-\left\langle A_nB_m,(S\varphi)(x+x_1,x-t+x_1-x_1^*)\right\rangle\\
&&
=-\left\langle\delta^{(n)}(x)\otimes\delta^{(m)}(t),(S\varphi)(x+x_1,x-t+x_1-x_1^*)\right\rangle\\
&&
=(-1)^{n+m+1}\d_x^n\d_t^m(S\varphi)(x+x_1,x-t+x_1-x_1^*)\Big|_{x=0,t=0}\\
&&
=\sum\limits_{j=0}^n \sum\limits_{i=0}^{n+m}F_{ji}(x,t)\d_x^j\d_t^i\varphi(x-t+x_1-x_1^*,x-t+x_1-x_1^*)\Big|_{x=0,t=0}\\
&&
=\sum\limits_{j=0}^n \sum\limits_{i=0}^{n+m}F_{ji}(0,0)\d_x^j\d_t^i\varphi(x_1,t_1^*)\\
&&
=\sum\limits_{j=0}^n \sum\limits_{i=0}^{n+m}F_{ji}(0,0)\left\langle\delta^{(j)}(x-x_1)\otimes\delta^{(i)}(t-t_1^*),
\varphi(x,t)\right\rangle.
\end{eqnarray*}
Here $F_{ji}(x,t)$ are known smooth functions of  $S$ and of all
its derivatives up to the order $n+m$.
Hence, for all $\psi(t)\in\D(0,T)$ we get
\begin{eqnarray*}
\lefteqn{
\left\langle\int\limits_0^L\delta^{(n)}(x-x_1)S(x,t)\delta^{(m)}(x-t-x_1^*)\,dx,\psi(t)\right\rangle}\\
&&
=\sum\limits_{j=0}^n \sum\limits_{i=0}^{n+m}F_{ji}(0,0)
\left\langle\int\limits_0^L\delta^{(j)}(x-x_1)\otimes\delta^{(i)}(t-t_1^*)\,dx,\psi(t)\right\rangle\\
&&
=\sum\limits_{j=0}^n \sum\limits_{i=0}^{n+m}F_{ji}(0,0)
\left\langle\delta^{(i)}(t-t_1^*)\otimes\delta^{(j)}(x-x_1),1(x)\otimes\psi(t)\right\rangle\\
&&
=\sum\limits_{i=0}^{n+m}F_{0i}(0,0)\left\langle\delta^{(i)}(t-t_1^*),\psi(t)\right\rangle.
\end{eqnarray*}
As a consequence,
\begin{eqnarray}
\label{eq:57_1}
\lefteqn{
I_0(t)=\int\limits_{t}^Lb_r(x)(S(x,t)a_r(x-t)+S_1(x,t))\,dx}\nonumber\\
&&
+(-1)^n(S(x,t)
a_r(x-t)
+S_1(x,t))_x^{(n)}|_{x=x_1}\\
&&
+(-1)^m
\d_x^mS(x+t+x_1^*,t)b_r(x+t+x_1^*)
\Big|_{x=0}
+\sum\limits_{i=0}^{n+m}F_{0i}(0,0)\delta^{(i)}(t-t_1^*).\nonumber
\end{eqnarray}
Observe that the first three summands in (\ref{eq:57_1}) are smooth
for $t>0$. Indeed, the second summand is smooth due to $a_r^{(i)}(0)=0$
for $0\le i\le n$ (see Assumption~1). The third summand is smooth due to
$b_r^{(i)}(L)=0$
for $0\le i\le m$ (see Assumption~2).

{\bf Further plan of the solution construction}.
We  split $\Omega_1^T$ into subdomains
$$
\Omega(i)=\{(x,t)\in\Omega_1^T\,|\,t-t_i^*<x<t-t_{i-1}^*\}
$$
and construct the solution separately in each $\Omega(i)$ and in a
neighborhood of each border between $\Omega(i)$ and $\Omega(i+1)$.
Here $t_0^*=0$, $1\le i\le k(T)$,
where $k(T)$  is defined
by
 inequalities $t_{k(T)}^*<T$ and
$t_{k(T)+1}^*\ge T$. The finiteness of $k(T)$  is obvious.

\section{Existence of the smooth solution on $\Omega(1)$}

\begin{lemma}\label{lemma:11}
There exists a  smooth solution to the problem
(\ref{eq:51})--(\ref{eq:53}) on $\Omega(1)$.
\end{lemma}
\begin{proof}
Under the assumption that $x_1^*<x_1$, we have $t_1^*<L$. Hence
$(x_1,t_1^*)\in\Omega_0$.
 Therefore any  solution which
is given by (\ref{eq:54}) on $\Omega_0^T$, 
is smooth on
$\Omega(1)$, and has the property given by Item 9 of Definition~\ref{defn:Omega}, satisfies the integral 
Volterra equation of the second kind
\begin{equation}\label{eq:57_2}
u(x,t)=S_3(x,t)+S_2(x,t)\int\limits_0^{t-x}b_r(\xi)u(\xi,t-x)\,d\xi,
\end{equation}
where
$$
S_2(x,t)=S(x,t)c_r(t-x)
$$
and
$$
S_3(x,t)=S_2(x,t)I_0(t-x)+S_1(x,t)
$$
are known by  (\ref{eq:57_1}).
The smoothness of $I_0(t-x)$  if $(x,t)\in\Omega(1)$ follows from
the facts that
$t-x<t_1^*$ and that $I_0(t)$ restricted to the
interval $(0,t_1^*)$ is smooth. Therefore $S_2$ and $S_3$ are smooth.

The lemma will follow from two claims. Set
$$
\Omega^{t(m)}(1)=\{(x,t)\in\Omega(1)\,|\,t<t(m)\}.
$$

{\it Claim 1. Given $m\in\N_0$, there exists a unique solution
$u\in C^m(\overline{\Omega^{t(m)}(1)})$ to the problem (\ref{eq:51})--(\ref{eq:53})
for some $t(m)>0$}. We apply the contraction principle to (\ref{eq:57_2}).
Comparing the difference of two continuous functions $u$ and
$\tilde u$ satisfying (\ref{eq:57_2}), we have
$$
|u-\tilde u|\le t(0)q
\max\limits_{(x,t)\in\overline{\Omega^{t(0)}(1)}}|u-\tilde u|,
$$
where
$$
q=\max\limits_{(x,t)\in\overline{\Omega(1)}}|S|
\max\limits_{t\in[0,t_1^*]}|c_r|
\max\limits_{x\in[0,L]}|b_r|.
$$
Choosing $t(0)<1/q$, we obtain the contraction property for the operator
defined by the right-hand side of (\ref{eq:57_2}).
The claim for $m=0$ follows.

Our next concern is the existence and uniqueness of a
$C^1(\overline{\Omega^{t(1)}(1)})$-solution for some $t(1)$. Let us consider
the problem
\begin{eqnarray}
\label{eq:ux}
\lefteqn{
\d_xu(x,t)=\d_xS_3(x,t)+\d_xS_2(x,t)\int\limits_0^{t-x}
b_r(\xi)u(\xi,\theta(x,t))\,d\xi}\nonumber\\
&&
-b_r(t-x)u(t-x,t-x)-S_2(x,t)\int\limits_0^{t-x}
b_r(\xi)(\d_tu)(\xi,t-x)\,d\xi.
\end{eqnarray}
From (\ref{eq:51}) we have $\d_tu=p(x,t)u+g(x,t)-\d_xu$.
We  choose an arbitrary $t(1)\le t(0)$.
Since $u$ is a known $C(\Omega^{t(1)}(1))$-function,
(\ref{eq:ux})  on $\overline{\Omega^{t(1)}(1)}$ is the Volterra integral equation of the second
kind with respect
to $\d_xu$. Assuming
in addition to the condition $t(1)\le t(0)$ that
$
t(1)<q,
$
we obtain the contraction property for (\ref{eq:ux}). On the account of
(\ref{eq:51}), the claim for $m=1$
follows.

Proceeding further by induction and using in parallel
(\ref{eq:51}), (\ref{eq:57_2}), and their suitable differentiations, we
complete the proof of the claim.

{\it Claim 2. In the domain
$\Omega^{t_1^*}(1)$ there exists a unique
smooth solution
 to the problem (\ref{eq:51})--(\ref{eq:53}).}
Given $m\in\N_0$, we prove that there exists a unique $u\in C^m(\Omega^{t_1^*}(1))$
in at most $\lceil t_1^*/t(m)\rceil$ steps by iterating the local existence
and uniqueness result in domains
$$
\Omega^{kt(m)}(1)\setminus\overline{\Omega^{(k-1)t(m)}(1)},\quad
1\le k\le\lceil T/t(m)\rceil.
$$
In particular, for $m=0$ in the $k$-th step of the proof  we have
\begin{equation}\label{eq:uuu}
\begin{array}{c}
\displaystyle
u(x,t)=S_3(x,t)+
S_2(x,t)\int\limits_0^{t-x-(k-1)t(m)}b_r(\xi)
u(\xi,t-x)\,d\xi\\
\displaystyle
+
S_2(x,t)\int\limits_{t-x-(k-1)t(m)}^{t-x}
b_r(\xi)
u(\xi,t-x)\,d\xi\\[7mm]
\displaystyle
\mbox{on}\,\,
\{(x,t)\in\Omega^{kt(m)}(1)\,|\,x\le t-(k-1)t(m)\}
\end{array}
\end{equation}
 and
\begin{equation}\label{eq:k}
u(x,t)=S(x,t)u(0,t-x)+S_1(x,t)\,\,\mbox{on}\,\,
\{(x,t)\in\Omega(1)\,|\,x\ge t-(k-1)t(m)\}.
\end{equation}
As in the latter formula $t-x\le(k-1)t(m)$,
the function $u$ defined by
(\ref{eq:k}) is smooth and
 known from the previous steps. This implies that
  the last summand in (\ref{eq:uuu})  is known and smooth. Hence
(\ref{eq:uuu}) is the Volterra integral equation of the second kind.
Applying now the argument used to prove Claim~1, we obtain
the existence and uniqueness of a continuous
solution $u$ to (\ref{eq:uuu}) on
$\Omega^{kt(m)}(1)\setminus\overline{\Omega^{(k-1)t(m)}(1)}$.
Since $k$ is an arbitrary integer in the range $1\le k\le\lceil T/t(m)\rceil$,
we have $u\in C\left(\Omega^{t_1^*}(1)\right)$. Further
we similarly proceed with all
derivatives f $u$.
Claim 2 is therewith proved.

The solution on the whole $\Omega(1)$ is now uniquely determined by the
formula
$$
u(x,t)=S(x,t)u(0,t-x)+S_1(x,t),
$$
where   $u(0,t-x)$ is a known smooth
function.
The latter is true due to
$0<t-x<t_1^*$ and Claim 2.

The proof of the claim is complete.
\end{proof}

From the
formulas (\ref{eq:54}) and
(\ref{eq:57_2}), Lemma~\ref{lemma:11}, and Assumption~1
it follows
 that $u$
is smooth in a  neighborhood
of the characteristic line
$x=t$.
This ensures that $u$ we  construct
satisfies Item 7 of Definition~\ref{defn:Omega}.

Under the assumption that $\Omega(2)$ is nonempty,
in the next section we give the formula of the solution on
$$
\Omega^{\varepsilon}(1)=\Omega(1)\cup
\{(x,t)\in\overline{\Omega(2)}\,|\,x>t-t_1^*-\varepsilon\}
$$
for a fixed $\varepsilon>0$ such that $t_1^*-\varepsilon>0$ and
\begin{equation}\label{eq:br2}
 b_r(x)=0,\quad
x\in[0,2\varepsilon].
\end{equation}
Such $\varepsilon$ exists by Assumption 2.

\section{The  solution on $\Omega^{\varepsilon}(1)$}

Write now
\begin{equation}\label{eq:v}
v(t)=\int\limits_{0}^{L}(b_r(x)+\delta^{(n)}(x-x_1))u\,dx=v_r(t)+v_s(t),
\end{equation}
where $v_r(t)$ and $v_s(t)$ are, respectively,
regular (smooth) and singular parts of $v(t)$. On the account of
(\ref{eq:u}), (\ref{eq:57_1}), (\ref{eq:br2}), and the fact that $x_1^*<x_1$,
we have on $[0,t_1^*+\varepsilon]$:
\begin{eqnarray}
\label{eq:vr}
\lefteqn{
v_r(t)=\int\limits_{t-t_1^*+\varepsilon}^{t}b_r(x)u(x,t)\,dx+
\int\limits_{t}^Lb_r(x)\left[S(x,t)a_r(x-t)+S_1(x,t)\right]\,dx}\nonumber\\
&&
+(-1)^n\d_x^n(S(x,t)a_r(x-t)
+S_1(x,t))|_{x=x_1}\nonumber\\
&&
+(-1)^m\d_x^m\left(S(x+t+x_1^*,t)b_r(x+t+x_1^*)\right)\Big|_{x=0}
\end{eqnarray}
and
\begin{equation}\label{eq:vs}
v_s(t)=\sum\limits_{i=0}^{n+m}F_{0i}(0,0)\delta^{(i)}(t-t_1^*).
\end{equation}
Note that the first summand in (\ref{eq:vr}) is a known smooth function.
This follows from the inclusion
$[t-t_1^*+\varepsilon,t]\times\{t\}\subset\Omega(1)
\cup\{(x,t)\,|\,x=t\}$, Lemma~\ref{lemma:11}
and Assumption 1. 

We distinguish two cases.

{\it Case 1. $t_1^*=t_1.$}
As easily seen from (\ref{eq:v}), (\ref{eq:vr}), and (\ref{eq:vs}),
$v(t)=v_r(t)$ on $[0,t_1^*+\varepsilon]$.
Thus, Item 6 of Definition~\ref{defn:Omega} for
$u$ we  construct is fulfilled. Furthermore,
\begin{equation}\label{eq:58}
u(0,t)=(\delta^{(j)}(t-t_1^*)+c_r(t))
v_r(t)=
\sum\limits_{i=0}^j
C_i\delta^{(i)}(t-t_1^*)+c_r(t)v_r(t)
\end{equation}
for $t\in(0,t_1^*+\varepsilon)$. The constants $C_i$  depend on
$v_r^{(k)}(t_1^*)$ for $0\le k\le j$ and can be computed by means of (\ref{eq:vr}).

{\it Case 2. $t_1^*\ne t_1.$}
Then $x_1-x_1^*=t_1^*$. Using
(\ref{eq:v}) and (\ref{eq:vr}),
we derive a similar formula for $u(0,t)$
on $(0,t_1^*+\varepsilon)$:
\begin{equation}\label{eq:59}
u(0,t)=c_r(t)\sum\limits_{i=0}^{n+m}F_{0i}(0,0)\delta^{(i)}(t-t_1^*)+c_r(t)v_r(t)
=\sum\limits_{i=0}^{n+m}E_{i}\delta^{(i)}(t-t_1^*)+c_r(t)v_r(t),
\end{equation}
where $E_i$ are constants depending on $F_{0,k}(0,0)$ and $c_r^{(k)}(t_1^*)$
for $0\le k\le n+m$.

Set
$$
Q(t)=\sum\limits_{i=0}^jC_i\delta^{(i)}(t-t_1^*)\quad\mbox{if}\,\,t_1^*= t_1
$$
 and
$$
Q(t)=\sum\limits_{i=0}^{n+m}E_{i}\delta^{(i)}(t-t_1^*)\quad\mbox{if}\,\,
t_1^*\ne t_1.
$$

\begin{lemma}\label{lemma:sol}
$u(x,t)$ given by the formula
\begin{equation}\label{eq:**}
u(x,t)=S(x,t)c_r(t-x)v_r(t-x)+S_1(x,t)+S(x,t)Q(t-x),
\end{equation}
where $v_r(t)$ is determined by (\ref{eq:vr}),
is
 a $\D'(\Omega)$-solution to the problem (\ref{eq:51})--(\ref{eq:53})
restricted to $\Omega^{\varepsilon}(1)$.
\end{lemma}
\begin{proof}
On the account of (\ref{eq:58}),
(\ref{eq:59}), and
the construction of the solution on $\Omega(1)$  done in
 Section 6,
it is enough to prove that the restriction of $S(x,t)Q(t-x)$
to
$Y=\{0\}\times(0,t_1^*+\varepsilon)$
is well defined and that $S(x,t)Q(t-x)$ satisfies (\ref{eq:51})
with $g(x,t)\equiv 0$ on $\Omega^{\varepsilon}
(1)$
in a distributional sense. The proof of the latter uses the argument as in the proof
of Lemma \ref{lemma:0}. To prove the former, consider the smooth bijective map
$$
\Phi:(x,t)\to (x,t-x-t_1^*).
$$
and its inverse
$$
\Phi^{-1}:(x,t)\to (x,x+t+t_1^*).
$$
 Applying Theorem 6, we have
$$
\WF(\Phi^*B_i)\subset \{(0,t+t_1^*,-\eta,
\eta),\eta\ne 0\}.
$$
Furthermore,
$$
N(Y)=\{(0,t,\xi,0)\}
$$
and therefore
$$
\WF(\Phi^*B_i)\cap N(Y)=\emptyset \quad
\mbox{for all} \quad 0\le i\le n+m.
$$
By Theorem~4, the restriction of $S(x,t)Q(\theta(x,t))$ to $Y$
is well defined. The lemma is therewith proved.
\end{proof}

\section{Construction of the
smooth solution on $\Omega(2)$}

To shorten notation, without loss of generality we assume that
$\max\overline{\Omega_1^T}\cap\{(x,t)\,|\,x=0\}\ge t_2^*$.
\begin{lemma}\label{lemma:12}
There exists a  smooth solution to the problem
(\ref{eq:51})--(\ref{eq:53}) on $\Omega(2)$.
\end{lemma}
\begin{proof}
 We start from the general formula of a smooth solution on $\Omega(2)$:
\begin{equation}\label{eq:59_1}
u(x,t)=S(x,t)u(0,t-x)+S_1(x,t).
\end{equation}
Since $S$ and $S_1$ are smooth, our task is  to prove that
there exists a smooth function identically equal to
$u(0,t-x)$
 on $\Omega(2)$. Since
$t_1^*<t-x<t_2^*$ if $(x,t)\in\Omega(2)$ and $c(t)=c_r(t)$ if $t\in(t_1^*,t_2^*)$,
it suffices to show the existence of a smooth function
$v_r(t)$ identically equal to $v(t)$
on $(t_1^*,t_2^*)$.
From the formula (\ref{eq:vs}) for $v_s(t)$ on $(0,t_1^*+\varepsilon)$
it follows that $v(t)=v_r(t)$
if $t\in(t_1^*,t_1^*+\varepsilon)$, where $\varepsilon$ is as in Section 7 and $v_r(t)$
is known and determined by
(\ref{eq:vr}).
To prove the lemma, it is sufficient to show
that there exists a smooth extension of $v_r(t)$ from $(0,t_1^*+\varepsilon)$
to $[t_1^*+\varepsilon,t_2^*)$ such that
$v_r(t)=v(t)$ if $t\in[t_1^*+\varepsilon,t_2^*)$.
If a such extension exists, then
by (\ref{eq:**})
it satisfies the following
integral equation
on $[t_1^*+\varepsilon,t_2^*)$:
\begin{equation}\label{eq:rv}
v_r(t)
=\int\limits_0^{t-t_1^*-\varepsilon}b_r(x)S(x,t)c_r(t-x)v_r(t-x)\,dx
+R(t),
\end{equation}
where
\begin{equation}\label{eq:R}
\begin{array}{cc}
\displaystyle
R(t)=\int\limits_{t-t_1^*-\varepsilon}^{P(t)}b_r(x)S(x,t)c_r(t-x)v_r(t-x)\,dx
+\int\limits_{0}^{P(t)}b_r(x)S_1(x,t)\,dx\\
\displaystyle
+I_0(t)+\int\limits_0^Lb_r(x)S(x,t)Q(t-x)\,dx,
\end{array}
\end{equation}
$$
P(t)=
\cases{t
&if
$L\le t$,
 \cr L
&if
$L\ge t
$, \cr}
$$
$b_r(x)$ is defined to be $0$ outside $[0,L]$, and
$v_r$ in the formula (\ref{eq:R}) is known and defined by (\ref{eq:vr}).
  One can easily see that the
first three summands in (\ref{eq:R}) are smooth functions on
$[t_1^*+\varepsilon,t_2^*)$. We now show that the last summand is a
$C^{\infty}[t_1^*+\varepsilon,t_2^*)$-function as well.
Indeed, take $\psi(t)\in\D(t_1^*+\varepsilon/2,t_1^*)$ and compute
\begin{eqnarray*}
\lefteqn{
\left\langle\int\limits_0^L b_r(x)S(x,t)\delta^{(j)}(t-x-t_1^*)\,dx,\psi(t)\right\rangle}\\
&&
=
\left\langle\delta^{(j)}(t-x-t_1^*),b_r(x)S(x,t)\psi(t)\right\rangle\\
&&
=-\left\langle\delta^{(j)}(x)\otimes 1(t),
b_r(t-x-t_1^*)S(t-x-t_1^*,t)\psi(t)\right\rangle\\
&&
=(-1)^{j+1}\left\langle\d_x^j\left(b_r(t-x-t_1^*)S(t-x-t_1^*,t)\right)
\Big|_{x=0},\psi(t)\right\rangle.
\end{eqnarray*}
We conclude that,
irrespective of whether
$t_1=t_1^*$
or $t_1\ne t_1^*$, the last summand
in (\ref{eq:R})  is a known  smooth function.
As follows from (\ref{eq:br2}),
the functions $v_r(t)$
defined by (\ref{eq:vr}) and
(\ref{eq:rv}) coincide at  $t=t_1^*+\varepsilon$.
The same is true
with respect to all the  derivatives of $v_r$.

Our task is therefore reduced to
show that there exists a  $C^{\infty}[t_1^*+\varepsilon,t_2^*)$-function
$v_r(t)$  satisfying (\ref{eq:rv}).
This follows from the fact that (\ref{eq:rv}) is the
integral Volterra equation of the second kind with respect to $v_r(t)$
(for  details see the proof of Lemma~\ref{lemma:11}).
The  proof is complete.
\end{proof}

\section{Completion of the construction}
Continuing our construction in this fashion,
we extend  $u$ over  a
neighborhood of each subsequent border between
$\Omega(i-1)$ and $\Omega(i)$
 and over
$\Omega(i)$ for all $3\le i\le k(T)$.
 Eventually
we construct $u$ on $\Omega^T$
for any $T>0$  in the sense of Definition~\ref{defn:Omega}
with $\Omega$ replaced by $\Omega^T$ and $\Pi$ replaced by
$\Pi^T=\{(x,t)\in\Pi\,|\,t<T\}$.
As easily seen from our construction, the
condition (\ref{eq:r}) is fulfilled with $\Omega_+$ and $\Omega_+'$
replaced by $\Omega^T\cap\Omega_+$ and $\Omega^T\cap\Omega_+'$,
respectively.
Since $T$ is arbitrary, the proof of Item 1 of Theorem~\ref{thm:exist}
is complete.
On the account of
Definition \ref{defn:Pi} and the definition of the restriction
$u\in\D'(\Omega)$ to a subset of $\Omega$ (see \cite[Section 5]{Horm}), Item 2
of Theorem~\ref{thm:exist} is a straightforward consequence of Item 1.
 Theorem~\ref{thm:exist} is therewith proved.

Assume that $S(x,t)\ne 0$ for all $(x,t)\in\Pi$. By
(\ref{eq:**}) it follows
from the  construction, that if
the singular part of $b(x)$ is the derivative of the 
Dirac measure of order
$n$,
then for each $i\ge 1$ there exist
$j>i$ and $n'\ge 1$  such that
$u$ is the derivative of the Dirac measure of order $n'$ along the
characteristic line $t-t_i^*$ and
$u$ is the derivative of the Dirac measure of order $n'+n$ along the
characteristic line $t-t_j^*$.
In contrast, this is not so if
singular parts of
the initial and the boundary data are Dirac measures.
In the latter case the solution preserves the same order of regularity
in time. Furthermore,
the assumption  $b_r^{(i)}(L)=0$ for all $i\in\N_0$ can be weakened to
$b_r(L)=0$.

Since $u$ restricted to $\Pi\setminus I_+$ is smooth,
Theorem~\ref{thm:order} follows from Item 2 of Theorem~\ref{thm:uniq}.

\section{Uniqueness of the  solution (Proof of
Theorem~\protect\ref{thm:uniq})}

The proof of Theorem~\protect\ref{thm:uniq} is based on 5 lemmas.

\begin{lemma}\label{lemma:uOmega0}
A $\D_+'(\Omega)$-solution $u$ to the problem (\ref{eq:51})--(\ref{eq:53})
is unique on $\Omega_0$.
\end{lemma}
\begin{proof}
Let $u$ and $\tilde u$ be two $\D_+'(\Omega_0)$-solutions
to the problem (\ref{eq:51})--(\ref{eq:52}). Then 
\begin{equation}\label{eq:L1}
\left \langle L(u-\tilde u),\varphi\right\rangle=\left\langle u-\tilde u,L^*\varphi\right\rangle=0
\quad \mbox{for all}\quad \varphi\in\D(\Omega_0),
\end{equation}
where
\begin{equation}\label{eq:LL}
  L=\d_t+\tilde{\lambda}\d_x-\tilde{p},\quad
L^*=-(\d_t+
\tilde{\lambda}\d_x+\tilde{p}).
\end{equation}
Our goal is to show that
\begin{equation}\label{eq:L5}
 \left\langle u-\tilde u,\psi\right\rangle=0\quad \mbox{for all}\quad
\psi\in\D(\Omega_0).
\end{equation}
Using the definition of
$\D_+'(\Omega_0)$ and (\ref{eq:L1}), it is sufficient to  prove that for every
$\psi\in\D(\Omega_0)$ there
exists $\varphi\in\D(\Omega_0)$ such that
\begin{equation}\label{eq:L4}
L^*\varphi=\psi\quad \mbox{on}\quad
\{(x,t)\in\Omega_0\,|\,t\ge 0\}.
\end{equation}
Fix $\psi\in\D(\Omega_0)$. If $\supp\psi\cap \{(x,t)\,|\,t>0\}=\emptyset$,
(\ref{eq:L5}) follows immediately from the definition of $\D_+'(\Omega_0)$.
We therefore assume that $\supp\psi\cap \{(x,t)\,|\,t>0\}\ne\emptyset$,
Consider the problem
$$
\varphi_t+\varphi_x=-p\varphi-\psi,\quad (x,t)\in
\{(x,t)\in\Omega_0\,|\,t>0\},
$$
$$
\varphi|_{t=0}=\varphi_0(x),\quad x\in(0,L),
$$
where $\varphi_0(x)\in\D(0,L)$ will be specified below.
This problem
has a unique smooth solution given by the formula
$$
\varphi(x,t)=\hat S(x,t)\varphi_0(x-t)+\hat S_1(x,t),
$$
where $\hat S_1$ is given by
(\ref{eq:S1}) with $p$ and $g$ replaced by $-p$ and $-\psi$,
respectively.

Fix $T(\psi)>0$ so that
$\supp\psi\cap\{(x,t)\,|\,t\ge T(\psi)\}=\emptyset$
and $\hat S(x,T(\psi))\ne 0$ for all $x$ with
$(x,T(\psi))\in\Omega_0$. The latter is ensured by
 Assumption 5.
Set
$$
\varphi_0(x-T(\psi))=
-\frac{\hat S_1(x,T(\psi))}{\hat S(x,T(\psi))}
$$
for $x$ such that $(x,T(\psi))\in\Omega_0$.
Changing coordinates $x\to\xi=x-T(\psi)$, we obtain
\begin{equation}\label{eq:phi0}
\varphi_0(\xi)=
-\frac{\hat S_1(\xi+T(\psi),T(\psi))}
{\hat S(\xi+T(\psi),T(\psi))}.
\end{equation}
We  construct the desired function $\varphi(x,t)$ by the formula
$$
\varphi(x,t)=
\cases{0
&if
$(x,t)\in\{(x,t)\in\Omega_0\,|\,t\ge T(\psi)\}$,
 \cr \hat S(x,t)\varphi_0(x-t)+\hat S_1(x,t)
&if
$(x,t)\in\{(x,t)\in
\Omega_0\,|\,0\le t\le T(\psi)\}$,
\cr \tilde\varphi(x,t)
&if
$(x,t)\in\{(x,t)\in
\Omega_0\,|\,t\le 0\}$,
 \cr}
$$
where $\tilde\varphi(x,t)$ is chosen so that $\varphi\in\D(\Omega_0)$.
The proof is complete.
\end{proof}

From now on we use a modified definition of $\Omega(i)$:
$$
\Omega(i)=\{(x,t)\in\Omega\,|\,t-t_i^*<x<t-t_{i-1}^*\},
\quad i\ge 1.
$$
Recall that $t_0^*=0$.

\begin{lemma}\label{lemma:uOmega1}
A $\D_+'(\Omega)$-solution
to the problem (\ref{eq:51})--(\ref{eq:53})
is unique on $\Omega(1)$.
 \end{lemma}
\begin{proof}
Assume that there exist two $\D_+'(\Omega)$-solutions
$u$ and $\tilde u$.
We will show that
\begin{equation}\label{eq:t2}
\left\langle v(t)-\tilde v(t),\psi(t)\right\rangle=0
\quad\mbox{for all}\quad
\psi(t)\in\D(0,t_1^*),
\end{equation}
where $v(t)$ is defined by Item 5 of Definition~\ref{defn:Omega} and
$\tilde v(t)$ is defined similarly
with $u$ replaced by $\tilde u$.  Postponing the proof, assume that
(\ref{eq:t2}) is true. Taking into account
Item 2 of Definition~\ref{defn:Omega0} and the fact that $c(t)=c_r(t)$
if $0<t<t_1^*$, we have
$$
 \left\langle L(u-\tilde u),\varphi\right\rangle=
\left\langle u-\tilde u,L^*\varphi\right\rangle=0\quad \mbox{for all}\quad \varphi\in\D(\Omega(1)).
$$
Let us prove that
\begin{equation}\label{eq:psi}
 \left\langle u-\tilde u,\psi\right\rangle=0\quad \mbox{for all}\quad
\psi\in\D(\Omega(1)).
\end{equation}
Following the argument used in the proof of Lemma~\ref{lemma:uOmega0},
 it is sufficient to show that,
given $\psi\in\D(\Omega(1))$, there exists $\varphi\in\D(\Omega(1))$
such that
$$
L^*\varphi=\psi\,\,\mbox{on}\,\,
\{(x,t)\in\Omega(1)\,|\,x\ge 0\}.
$$
We concentrate on the case that $\supp\psi\cap\{(x,t)\in\Omega(1)\,|\,x>0\}
\ne\emptyset$.
Otherwise (\ref{eq:psi}) is immediate
because $u-\tilde u\in\D_+'(\Omega(1))$.
Consider the problem
$$
\varphi_t+\varphi_x=-p\varphi-\psi,\quad (x,t)\in
\{(x,t)\in\Omega(1)\,|\,x>0\},
$$
$$
\varphi|_{x=0}=\varphi_1(t),\quad t\in(0,t_1^*),
$$
where $\varphi_1(t)\in\D(0,t_1^*)$ is a fixed function. Let
$T(\psi)>0$ be the same as in the proof of Lemma~\ref{lemma:uOmega0}.
We specify $\varphi_1(\xi)$ by
\begin{equation}\label{eq:phi1}
\varphi_1(\xi)=
-\frac{\hat S_1(T(\psi)-\xi,T(\psi))}
{\hat S(T(\psi)-\xi,T(\psi))}
\end{equation}
and construct the desired $\varphi$ similarly to the construction of $\varphi$
in the proof of Lemma~\ref{lemma:uOmega0}.
To finish the proof of the lemma, it remains to show that
\begin{equation}\label{eq:36_0}
\left\langle v-\tilde v,\psi(t)\right\rangle=0
\quad\mbox{for all}\quad
\psi(t)\in\D(\varepsilon i,\varepsilon i+2\varepsilon),
\end{equation}
for each $0\le i\le t_1^*/\varepsilon-2$, where $\varepsilon>0$ is chosen so that
$t_1^*/\varepsilon$ is an integer and
\begin{equation}\label{eq:br}
b_r(x)=0\,\,\mbox{for}\,\,x\in[0,2\varepsilon].
\end{equation}
Such $\varepsilon$ exists by  Assumption 2. We prove (\ref{eq:36_0}) by induction on $i$.

{\it Claim 1 (the base case). (\ref{eq:36_0}) is true for
$i=0$.} We will use the following representations for $u$ and $\tilde u$
on $\Omega_+$ which  are possible owing to  Item 3 of
Definition~\ref{defn:D+}:
\begin{equation}\label{eq:u01}
\begin{array}{c}
u=u_0+u_1\,\,\mbox{in}\,\,\D'(\Omega_+),\\\displaystyle
\tilde u=\tilde u_0+\tilde u_1\,\,\mbox{in}\,\,\D'(\Omega_+),
\end{array}
\end{equation}
where $u_0=u$ and $\tilde u_0=\tilde u$ in $\D'(\Omega_0\cap\Omega_+)$,
$u_0=\tilde u_0\equiv 0$ on $\overline{\Omega_1}\cap\Omega_+$, $u_1=u$ and
$\tilde u_1=\tilde u$ in $\D'(\Omega_1\cap\Omega_+)$,
$u_1=\tilde u_1\equiv 0$  on $\overline{\Omega_0}\cap\Omega_+$.

We first prove that
\begin{equation}\label{eq:t3}
\left\langle v-\tilde v,\psi(t)\right\rangle=\left\langle u_1-\tilde u_1,b_r(x)\psi(t)\right\rangle
\quad\mbox{for all}\quad
\psi(t)\in\D(0,4\varepsilon).
\end{equation}
Accordingly to  Item 1 of Definition~\ref{defn:D+},
\begin{equation}\label{eq:uu}
\begin{array}{c}
\left\langle v-\tilde v,\psi(t)\right\rangle=\left\langle u-\tilde ub(x),1(x)\otimes\psi(t)\right\rangle
\\\displaystyle
=\left\langle (u_0-\tilde u_0)b(x),1(x)\otimes\psi(t)\right\rangle
+\left\langle (u_1-\tilde u_1)b(x),1(x)\otimes\psi(t)\right\rangle,
\end{array}
\end{equation}
where $b_r(x)=0$, $x\not\in[0,L]$. By
Lemma~\ref{lemma:uOmega0}, $u_0=\tilde u_0$ in $\D'(\Omega_0\cap\Omega_+)$.
Applying in addition Item 1 of Theorem~\ref{thm:exist} and
Proposition~\ref{prop:p},
we have
\begin{equation}\label{eq:I}
\left\langle(u_0-\tilde u_0)(x,t)b(x),1(x)\otimes\psi(t)\right\rangle=
\left\langle I_0(t)-\tilde I_0(t),\psi(t)\right\rangle,
\end{equation}
where $I_0(t)$ is defined by
(\ref{eq:I_0}) and $\tilde I_0(t)$ is defined by
(\ref{eq:I_0})  with $u_0$ replaced by~$\tilde u_0$.
From (\ref{eq:57_1}) we have $I_0(t)=\tilde I_0(t)$ for
$0<t<4\varepsilon$. Hence the right-hand side of (\ref{eq:I}) is equal to 0.
On the  account of the inclusions $\supp(u_1-\tilde u_1)\subset\overline{\Omega_1}$ and
$\supp\psi(t)\subset[0,2\varepsilon]$, (\ref{eq:uu}) does not depend
on $b(x)$ outside  $[0,2\varepsilon]$.
Since  $x_1^*<x_1$, $b(x)=b_r(x)$ on $[0,2\varepsilon]$. Therefore
(\ref{eq:uu}) implies (\ref{eq:t3}).   Claim 1 now
follows from (\ref{eq:br}).

Assume that (\ref{eq:36_0}) is true for
$i=k-1$, $k\ge 1$  and prove that it
is true for $i=k$.

{\it Claim 2 (the  induction step). (\ref{eq:36_0}) is true for
$i=k$, $k\ge 1$.}
The proof is similar to the proof of Claim 1.
Based on the induction assumption and applying the argument used in the
proof of (\ref{eq:psi}), we obtain
\begin{equation}\label{eq:eq}
u=\tilde u\,\,\mbox{in}\,\,\D_+'(G(k-1)),
\end{equation}
where
$$
G(k)=\Omega(1)\cap\{(x,t)\,|\,x>t-\varepsilon k-2\varepsilon)\}.
$$
Applying in addition Item 1 of Theorem~\ref{thm:exist},
Proposition~\ref{prop:p}, and Lemma~\ref{lemma:11}, we conclude that $u$
is smooth on $G(k-1)\cap\Omega_+$.
Owing to (\ref{eq:eq}) and the latter fact,
 the following  representations for
$u$ and $\tilde u$ on $\Omega_+$ are possible:
$$
u=u_0+u_{k-1}+u_k\,\,\mbox{in}\,\,\D'(\Omega_+),
$$
$$
\tilde u=u_0+u_{k-1}+\tilde u_k\,\,\mbox{in}\,\,\D'(\Omega_+),
$$
where $u_0$ is the same as in
(\ref{eq:u01}),
$u_{k-1}=u$  in
$\D'(G(k-1)\cap\Omega_+)$,
$u_{k-1}\equiv 0$  on
$\Omega_+\setminus G(k-1)$, $u_k=u$ and
$\tilde u_k=\tilde u$ in
$\D'(\Omega_+\setminus(\overline{G(k-1)}\cup\overline{\Omega_0}))$,
$u_k=\tilde u_k\equiv 0$  on
$\Omega_+\cap(\overline{G(k-1)}\cup\overline{\Omega_0})$.
Similarly to (\ref{eq:t3}), we derive the equality
$$
\left\langle v-\tilde v,\psi(t)\right\rangle=\left\langle u_k-\tilde u_k,b_r(x)\psi(t)\right\rangle
\quad\mbox{for all}\quad
\psi(t)\in\D(\varepsilon k,\varepsilon k+2\varepsilon).
$$
The claim follows from the
 support properties of $u_k-\tilde u_k$, $\psi(t)$,
and $b_r$ given by~(\ref{eq:br}).

The proof is complete.
\end{proof}

Set
$$
\Omega^{\varepsilon}(0,1)=
\{(x,t)\in\Omega\,|\,x-\varepsilon<t<x+\varepsilon)\}.
$$

\begin{lemma}\label{lemma:uOmega0Omega1}
A $\D_+'(\Omega)$-solution
to the problem (\ref{eq:51})--(\ref{eq:53})
is unique on $\Omega^{\varepsilon}(0,1)$ provided $\varepsilon$ is small enough.
 \end{lemma}
\begin{proof}
Let $u$ and $\tilde u$ be two
$\D_+'(\Omega)$-solutions
to the problem
(\ref{eq:51})--(\ref{eq:53}).
Fix $\varepsilon>0$ so that the condition (\ref{eq:br}) is fulfilled.
By Claim 1 in the proof of Lemma~\ref{lemma:uOmega1},
(\ref{eq:36_0}) is true for $i=0$. Therefore
$$
\left\langle L(u-\tilde u),\varphi\right\rangle=\left\langle u-\tilde u,L^*\varphi\right\rangle=0\quad
\mbox{for all}\quad
\varphi\in\D(\Omega^{\varepsilon}(0,1)).
$$
Our task is to prove (\ref{eq:psi}) with $\Omega(1)$ replaced by
$\Omega^{\varepsilon}(0,1)$. In fact, we prove that,
given $\psi\in\D(\Omega^{\varepsilon}(0,1))$,
there exists $\varphi\in\D(\Omega^{\varepsilon}(0,1))$ satisfying the initial
boundary  problem
$$
\varphi_t+\varphi_x=-p\varphi-\psi,\quad
(x,t)\in\Omega^{\varepsilon}(0,1)\cap\Omega_+,
$$
$$
\varphi|_{t=0}=\varphi_0(x),\quad x\in[0,\varepsilon),
$$
$$
\varphi|_{x=0}=\varphi_1(t),\quad t\in[0,\varepsilon).
$$
Here $\varphi_0(x)\in C^{\infty}[0,\varepsilon)$  is a fixed function identically
equal to 0 in a neighborhood of $\varepsilon$,
$\varphi_1(t)\in C^{\infty}[0,\varepsilon)$  is a fixed function identically
equal to 0 in a neighborhood of
$\varepsilon$, and $\varphi_0^{(i)}(0)=\varphi_1^{(i)}(0)$ for all $i\in\N_0$.
We construct $\varphi(x,t)$, combining the constructions of $\varphi(x,t)$
in the proofs of Lemmas~\ref{lemma:uOmega0} and~\ref{lemma:uOmega1}.
Thus we fix
$T(\psi)>0$ to be the same as in the proof of Lemma~\ref{lemma:uOmega0}
and specify $\varphi_0(x)$ and $\varphi_1(t)$  by
(\ref{eq:phi0}) and (\ref{eq:phi1}), respectively.
Let
$$
\varphi(x,t)
$$
$$
=\cases{0
&if
$(x,t)\in\{(x,t)\in\Omega^{\varepsilon}(0,1)\,|\,t\ge T(\psi)\}$,
 \cr \hat S(x,t)\varphi_0(x-t)+\hat S_1(x,t)
&if
$(x,t)\in\{(x,t)\in
\overline{\Omega_0}\cap\Omega^{\varepsilon}(0,1)\,|\,0\le t\le T(\psi)\}$,
 \cr \hat S(x,t)\varphi_1(t-x)+\hat S_1(x,t)
&if
$(x,t)\in\{(x,t)\in
\overline{\Omega(1)}\cap\Omega^{\varepsilon}(0,1)\,|\,0\le t\le T(\psi)\}$,
\cr \tilde\varphi(x,t)
&if
$(x,t)\in\{(x,t)\in
\Omega^{\varepsilon}(0,1)\,|\,x\le 0\,\mbox{or}\,t\le 0\}$,
 \cr}
$$
where $\tilde\varphi(x,t)$ is chosen so that $\varphi\in\D(\Omega^{\varepsilon}(0,1))$.

The proof is complete.
\end{proof}

For every $i\ge 1$ fix $\varepsilon_i$  such that
 $t_i^*-\varepsilon_i>t_{i-1}^*$, $t_i^*+\varepsilon_i<t_{i+1}^*$, and
\begin{equation}\label{eq:bri}
b_r(x)=0\,\,\mbox{for}\,\,x\in[0,4\varepsilon_i].
\end{equation}
Set
$$
Q(i)=\{(x,t)\,|\,t-t_i^*-\varepsilon_i<x<t-t_i^*+\varepsilon_i\}.
$$

\begin{lemma}\label{lemma:uOmega12}
A $\D_+'(\Omega)$-solution
to the problem (\ref{eq:51})--(\ref{eq:53})
is unique on $Q(1)$.
\end{lemma}
\begin{proof}
Assume that there exist two $\D_+'(\Omega)$-solutions
$u$ and $\tilde u$
and show that
\begin{equation}\label{eq:L25}
\left\langle v-\tilde v,\psi(t)\right\rangle=0\quad\mbox{for all}\quad
\psi(t)\in\D(t_1^*-\varepsilon_1,t_1^*+\varepsilon_1).
\end{equation}
By Lemmas~\ref{lemma:11}
and~\ref{lemma:uOmega1}, Item 1 of Theorem~\ref{thm:exist}, and
Proposition~\ref{prop:p}, any solution to      (\ref{eq:51})--(\ref{eq:53})
restricted to $\Omega(1)$ is smooth. Based on this fact and on
Lemmas~\ref{lemma:uOmega0}--\ref{lemma:uOmega0Omega1},
similarly to (\ref{eq:t3}), we derive the
equality
$$
\left\langle v-\tilde v,\psi(t)\right\rangle=
\left\langle u_1-\tilde u_1,b_r(x)\psi(t)\right\rangle
$$
$$
\mbox{for all}\,\,
\psi(t)\in\D(t_1^*-\varepsilon_1,t_1^*+\varepsilon_1),
$$
where $u_1=u$ and $\tilde u_1=
\tilde u$ in $\D'(G)$,  $u_1$ and $\tilde u_1$ are identically
equal to zero on $\Omega_+\setminus G$. Here
$$
G=\{(x,t)\in\Omega_+\,|\,x<t-t_1^*+\varepsilon_1)\}).
$$
The equality (\ref{eq:L25}) now
follows
from the support properties of $u_1-\tilde u_1$, $\psi$, and $b_r$
given by (\ref{eq:bri}) for $i=1$.

We further distinguish two cases.

{\it Case 1. $t_1^*\ne t_1$}. Then $c(t)=c_r(t)$ for $t$ in the range
$t_1^*-\varepsilon_1<t<t_1^*+\varepsilon_1$. Applying
(\ref{eq:L25}) and Item 2 of Definition~\ref{defn:D+}, we have
\begin{equation}\label{eq:Lu}
L(u-\tilde u)=0\,\,\mbox{in}\,\,
\D'(Q(1)).
\end{equation}

{\it Case 2. $t_1^*=t_1$}.
Then $c(t)=\delta^{(j)}(t-t_1)+c_r(t)$. By
Item 6 of Definition~\ref{defn:Omega},
$v-\tilde v$ is smooth in a neighborhood of $t_1^*$.
Combining the latter with (\ref{eq:L25}), we get (\ref{eq:Lu}).

In the rest of the proof we proceed as in  the proof of
 Lemma~\ref{lemma:uOmega1}.
\end{proof}

\begin{lemma}\label{lemma:uOmega2}
A $\D_+'(\Omega)$-solution
to the problem (\ref{eq:51})--(\ref{eq:53})
is unique on $\Omega(2)$.
\end{lemma}
\begin{proof}
We follow the proof of Lemma~\ref{lemma:uOmega1} with $\Omega(1)$ replaced by
$\Omega(2)$ and with minor changes caused by the fact that due to
 Lemmas~\ref{lemma:uOmega12} and~\ref{lemma:12},
$u$ and $\tilde u$ are smooth on
$\Omega(2)\cap\Omega_+\cap\{(x,t)\,|\,
x>t-t_1^*-\varepsilon_1\}$. Hence
(\ref{eq:t2}) is true with $(0,t_1^*)$ replaced by $(t_1^*+\varepsilon_1/2,t_2^*)$.
\end{proof}

Continuing in this fashion, we eventually prove the uniqueness over subsequent
$\Omega(i)$ and $Q(i)$ for any desired  $i\in\N$.
Summarizing it with Lemmas~\ref{lemma:uOmega0} and~\ref{lemma:uOmega0Omega1} and
Theorem~\ref{thm:shilow}, we obtain Item 1 of Theorem~\ref{thm:uniq}.

Item 2 of Theorem~\ref{thm:uniq}
is a straightforward consequence of
Item~1 of Theorem~\ref{thm:uniq}, Item 2 of Theorem~\ref{thm:exist},
and  Proposition~\ref{prop:p}.

\subsection*{Acknowledgments}
The author is  thankful to the members of the DIANA group  for their kind
hospitality during her stay at the Vienna university.

\qquad
\parbox{60mm}{
\small I. Kmit \\
Institute for Applied Problems\\
 of Mechanics and Mathematics,
\\
Ukrainian Academy of Sciences\\
Naukova St.\ 3b,\\
 79060 Lviv, Ukraine
\\
E-mail: kmit@informatik.hu-berlin.de
}

\end{document}